\documentclass[11pt]{amsart}
\usepackage[a4paper, left=2.5cm, right=2.5cm, top=3cm, bottom=3cm]{geometry}

\usepackage{latexsym,amssymb,amsmath,graphics}
\usepackage[dvips]{graphicx}
\usepackage{color}
\usepackage{hyperref}

\definecolor{orange}{rgb}{1.0, 0.5, 0.0}

\definecolor{uuuuuu}{rgb}{0.26666666666666666,0.26666666666666666,0.26666666666666666}

\theoremstyle{plain}
\newtheorem{thm}{Theorem}
\newcounter{counter}

\newtheorem{prop}{Proposition}

\newtheorem{rmk}[prop]{Remark}

\newcommand {\R} {\mathbb{R}} \newcommand {\Z} {\mathbb{Z}}
 
\newcommand {\C} {\mathbb{C}}

\numberwithin{equation}{section}

\allowdisplaybreaks

\title{A theorem concerning Fourier transforms: a survey}

\date{\today}

\author[A. Fern\'andez-Bertolin]{Aingeru Fern\'andez-Bertolin}
\address[Aingeru Fern\'andez-Bertolin]{Universidad del Pa\'is Vasco / Euskal Herriko Unibertsitatea,
	48080 Bilbao, Spain}
\email{aingeru.fernandez@ehu.eus}
\author[L. Vega]{Luis Vega}
\address[Luis Vega]{BCAM -- Basque Center for Applied Mathematics,
	48009 Bilbao, Spain,
	\newline \phantom{\quad} \&
	Universidad del Pa\'is Vasco / Euskal Herriko Unibertsitatea,
	48080 Bilbao, Spain}
\email{luis.vega@ehu.eus }



\begin{document}

\maketitle

\begin{abstract}
In this note, we highlight the impact of the paper G.\ H.\ Hardy, {\em A theorem concerning Fourier transforms}, J. Lond. Math. Soc. (1) {\bf 8} (1933), 227--231
in the community of harmonic analysis in the last 90 years, reviewing, on the one hand, the direct generalizations of the main results and, on the other, the different connections to related areas and new perspectives.
\end{abstract}

\section{Introduction}

A century ago, in 1925, the mathematical physics community enjoyed a breakthrough that since then has had unlimited consequences: Heisenberg showed that non-commutativity is the key to quantum theory. Soon after that, in 1927, the Heisenberg uncertainty principle was established in \cite{He27}, exploiting the fact that the commutator of the position and momentum observables is the identity.

Mathematically speaking, we can state Heisenberg uncertainty principle through the so-called Heisenberg--Pauli--Weyl inequality by Kennard and Weyl, \cite{Ke27,We28}:
\begin{equation}\label{hpwineq}
\int_{\R^d}|x|^2|f(x)|^2\int_{\R^d}|\xi|^2|\widehat{f}(\xi)|^2\ge
 \frac{d^2}{4}\|f\|_2^4,\ 
 \forall f\in L^2(\R^d),
\end{equation}
where $\|\cdot\|_2$ stands for the $L^2(\R^d)$ norm. It is not difficult to check that the only extremizers of this inequality are Gaussians.

This last inequality reflects the fact, pointed out by Norbert Wiener (see \cite{Ha33,Th}), that a function and its Fourier transform, which we define in this paper as
\[
\widehat{f}(\xi)=\frac{1}{(2\pi)^{d/2}}\int_{\R^d}f(x)e^{-ix\cdot\xi}dx,
\]
cannot be simultaneously very localized. In other words, the more concentrated in space, the more dispersed in frequency. This vague statement is what we call the \textit{uncertainty principle}, as different interpretations of localization lead to different inequalities sharing this philosophy. We would like to point out the nice survey \cite{WW21} where different uncertainty principles are stated and linked to a general and somewhat simple inequality, called by the authors the \textit{primary uncertainty principle}.

The goal of this note is to emphasize the mathematical impact that one of the interpretations of the uncertainty principle, made by G. H. Hardy not so long after Heisenberg's achievement, has had over the last 100 years. In 1933, Hardy shows in the paper \textit{A theorem concerning Fourier transforms}, \cite{Ha33}, the following one-dimensional results:

\begin{thm}\label{generalhup}
    If $f$ and $\widehat{f}$ are both $O(|x|^m e^{-x^2/2})$ for large $x$ and some $m$, then each is a finite liner combination of Hermite functions.
\end{thm}

\begin{thm}\label{classicalhup}
In particular, if $f$ and $\widehat{f}$ are both $O(e^{-x^2/2})$ then $f=\widehat{f}=Ae^{-x^2/2}$, where $A$ is a constant; and if one is $o(e^{-x^2/2})$, then both are null. 
\end{thm}

\begin{rmk}\label{hardyab}
These statements are written as in the original manuscript of Hardy. It is easy to check from Theorem \ref{classicalhup} that if $f$ is $O(e^{-ax^2})$, $\widehat{f}$ is $O(e^{-bx^2})$, and $ab>1/4$ then both are null, having the characterization $f(x)=C e^{-ax^2}$ in the case $ab=1/4$. Usually this is what we call \textit{Hardy's uncertainty principle}.
\end{rmk}

Using Hardy's words, this gives \textit{the most precise interpretation possible of Wiener's remark.}

%It is worth pointing out that the value $m$ in Theorem \ref{generalhup} is related to the degree of Hermite functions in the linear combination. Furthermore, Hermite functions are also connected to the Heisenberg--Pauli--Weyl inequality: It is a nice exercise to check that if we restrict the $L^2(\R^d)$ space to the subspace of functions orthogonal to the first $k$ Hermite functions, then the best constant achieved in \eqref{hpwineq} is improved, quantifying this improvement with respect to $k$. Also, the extremizer in this space is given by a Hermite function.

It is not the purpose of this note to give details of the proof of these results. In his short-but-brilliant paper, Hardy provides two different proofs. The first one nowadays appears in different textbooks, sometimes as an exercise with hints to help the reader, see for instance \cite{DM,HJ,SS,Th}. It is based on the Phragm\'en-Lindel\"of principle and therefore on a log-convexity argument. This convexity argument will play a relevant role in Section 4 below, where a more general version of Phragm\'en-Lindel\"of due to Agmon and Nirenberg \cite {Ag66} \footnote{ We want to thank N. Garofalo for calling our attention to this work.}is used. 

The second approach, more elementary,   apparently did not catch the attention of the mathematical community and was forgotten. We refer the reader to the survey \cite{FBM21} of the first author and Malinnikova, where an argument in the flavor of Hardy's second proof, also inspired by the work of Hedenmalm \cite{He12}, is used to give an alternative proof of a uniqueness result for Schr\"odinger evolutions equivalent to Hardy's uncertainty principle. We will review some of the generalizations of Hardy's uncertainty principle through the last 90 years, to see how this result has attracted the interest of several researchers and is a principle that, in some sense, is still alive in the literature and connected to hot topics in mathematics. We note that there are already many surveys in the literature on uncertainty principles, including the one by Hardy, \cite{BD06, FBM21, FS97, Ja06}. In some sense, our note can be seen as an update of the previous ones. The interested reader should also read the book by Thangavelu \cite{Th}, focused on Hardy's theorem. Due to the large amount of results and research directions related to Hardy's uncertainty principle, we will very briefly describe some of these directions.

The rest of the note is organized as follows: In Section 2 we will derive Heisenberg's uncertainty principle and other well-known inequalities from an abstract uncertainty principle.
In Section 3 we study generalizations and extensions that are stated under the same framework of the original result. In Section 4 we comment on the search of a real variable approach towards Hardy's uncertainty principle, that links the classical formulation to a dynamical version of the principle. Finally, in Section 5 we discuss some other extensions of Hardy's uncertainty principle to settings that differ from $\R$ or $\R^d$.

\textbf{Acknowledgments}. A. Fern\'andez-Bertolin is supported by the project IT1615-22 (Basque Government) and by the project PID2021-122156NB-I00 funded by\newline MICIU/AEI/10.13039/501100011033.

 L. Vega is supported by MICIN (Spain) CEX2021-001142, by PID2021-126813NB-I00 funded by MICIU/AEI/10.13039/501100011033 and by "ERDF A way of making Europe”, and IT1615-22 (Basque Government).

We thank the anonymous referees for several remarks that improved the presentation of the manuscript.

\section{An abstract Uncertainty Principle}

The idea behind the uncertainty princiciple is that the multiplication of matrices or in general the composition of linear operators is  non-commutative. That is to say that the commutator  associated to two linear operators is typically non-zero and as a consequence an inequality can be proved. 

In a completely formal way one can proceed as follows. 
Consider two operators,
one, $ \mathcal{S}$ symmetric
$$ \langle \mathcal{S}f,f\rangle= \langle f,\mathcal{S}f\rangle,$$
and the other one,
$\mathcal{A\,} $ skewsymmetric
$$\langle \mathcal{A}f,f\rangle =  -\langle f,\mathcal{A}f\rangle,$$ 
so that,
$$\langle (\mathcal S+\mathcal A)f,(\mathcal S+\mathcal A)f\rangle=\langle \mathcal S f,\mathcal S f\rangle+\langle \mathcal A f,\mathcal A f\rangle+
\langle (\mathcal S \mathcal A-\mathcal A 	\mathcal S)f,f\rangle\ge 0.$$

Hence
$$\langle (
\mathcal A\mathcal S-\mathcal S\mathcal A)f,f\rangle\le \|\mathcal S f\|_{L^2}^2+\|\mathcal A f\|_{L^2}^2.$$

Now adjusting the inequality for any real $\lambda$, changing $\mathcal A$ into $ \pm\lambda \mathcal A$ and $\mathcal S$ into $\frac{1}{\lambda} \mathcal S$, and then choosing appropriately $\lambda$ we get the inequality
$$\left|\langle (
\mathcal A\mathcal S-\mathcal S\mathcal A)f,f\rangle\right|\le 2\|\mathcal S f\|_{L^2}\|\mathcal A f\|_{L^2}.$$

A  further generalization can be obtained just by taking $\ A-\langle \mathcal A f,f\rangle 1$ and $S-\langle \mathcal S f,f\rangle 1$. 

The above abstract inequality turns out to be very useful as the next three examples illustrate. The justification of the different steps needed for the corresponding proofs is nowadays quite standard and we omit it.

{\bf Example 1.}\\
Choose in one dimension $\mathcal S f=xf$ and $\mathcal A f=f^\prime$, so that
$\mathcal A\mathcal S - \mathcal S\mathcal A=\displaystyle\frac d{dx}x-x\frac d{dx}=1$. The consequence is the well know Heisenberg uncertainty principle
$$\|f\|^2_{L^2}\le 2 \|f^\prime\|_{L^2}\,\|xf\|_{L^2}.$$

The extremizer is easily obtained:
$$(\mathcal S+\mathcal A)f=0\quad\Longleftrightarrow\quad f(x)=c e^{-x^2/2}.$$

Moreover, the extremizer  satisfies
$$(\mathcal S-\mathcal A)(\mathcal S+\mathcal A)f=
-f^{\prime\prime}+x^2f-f=0,$$
and therefore can be also characterised as the eigenvector associated to the first eigenvalue of the (quantum) harmonic oscillator. Finally applying iteratively $\mathcal S-\mathcal A= x-\frac{d}{dx} $ to $f$ all the rest of eigenvalues are easily obtained and the corresponding eigenfunctions can be written in terms of Hermite functions. This procedure generalizes without any difficulty to higher dimensions. The above argument can be used  to check that if we restrict the $L^2(\R^d)$ space to the subspace of functions orthogonal to the first $k$ Hermite functions, then the best constant achieved in \eqref{hpwineq} is improved depending on $k$ and the extremizer in this space is given by a Hermite function. We emphasize that Hermite functions also appear naturally in Hardy's uncertainty principle, Theorem \ref{generalhup}. The value $m$ in Theorem \ref{generalhup} is related to the degree of Hermite functions in the linear combination.

{\bf Example 2.}\\ In this case we consider the dimension $d\geq2$,  $\mathcal S f=\displaystyle\frac x{|x|}$, and  $\mathcal A f=\nabla f$. Then
$$\mathcal A\mathcal S -\mathcal S\mathcal A=\displaystyle\frac {d-1}{|x|}$$
and
$$\displaystyle\int\frac{d-1}{|x|}|f|^2\le 2\left(\int|f|^2\right)^{1/2}
\left(\int |\nabla f|^2\right)^{1/2}.$$
The extremizer can be also explicitly computed in this case 
$$(\mathcal S+\mathcal A)f=0\quad\Longleftrightarrow\quad f(x)=ce^{-|x|},
$$
and appears as the eigenvector of the first eigenvalue of the Schr\"odinger equation for the Coulombian potential,
$$(\mathcal S-\mathcal A)(\mathcal S+\mathcal A)f=\Delta f-\displaystyle\frac{d-1}{|x|}f+f=0.$$

{\bf Example 3.}\\ Take $d\geq3$, $\mathcal S f=\displaystyle\frac {d-2}2\frac x{|x|^2}f$, and $\mathcal A f=\nabla f$. Then
$$\mathcal A\mathcal S -\mathcal S\mathcal A=\displaystyle\frac {(d-2)^2}{2}\frac 1{|x|^2}f,$$
and we get the well known Hardy’s Inequality 
$$\left(\displaystyle\int\frac{(d-2)^2}{2|x|^2}|f|^2\right)^{1/2}\le
\left(2\int\left|\nabla f\right|^2\right)^{1/2}.$$
Also,
$$(\mathcal S+\mathcal A)f=0\quad\Longleftrightarrow\quad f=|x|^{1-\frac d2},$$
so that, since $\nabla f\notin L^2$, the best constant of the inequality is not achieved. Finally,
$$(\mathcal S-\mathcal A)(\mathcal S+\mathcal A)f=-\Delta f-\displaystyle\frac{(d-2)^2}{4|x|^2}f=0.$$
It is well known that the constant $\frac{(d-2)^2}{4}$ is a natural threshold for self-adjointness of Schr\"odinger operators with potentials that have a point-wise quadratic singularity like the one above.

\section{Generalizations of Hardy's uncertainty principle}
Soon after Hardy's result, Morgan generalized Theorem \ref{classicalhup}, stating what we could call Morgan's uncertainty principle, \cite{Mo34}:

\begin{thm}\label{morgan}
Let $p>2$ and $A>0$. If $f$ is $O(e^{-A x^p})$ and $\widehat{f}$ is $O(e^{-B x^{q}})$, with $q$ the conjugate exponent of $p$ and $B>\frac{\sin(\pi/(2(p-1))}{q (Ap)^{q-1}}$, then both are null.
\end{thm}

One can check from Morgan's research the optimality of this result. Also, note that the case $p=2,\ A=1/2$ leads to the condition $\widehat{f}=O(-B x^2)$ with $B>1/2$. In other words, this result generalizes the one by Hardy.

Several years after Hardy's and Morgan's work, several extensions were made. We point here to three different directions:

\begin{itemize}
    \item Cowling and Price state $L^p$ versions of Theorem \ref{classicalhup}, still in $\mathbb{R}$, \cite{CP83}. It is worth mentioning that out of the $L^\infty$ setting, $f$ and $\widehat{f}$ are both null in the case $ab\ge 1/4$, as opposed to the $L^\infty$ version where $ab=1/4$ implies that $f$ is a multiple of the Gaussian $e^{-ax^2}$.
    \item In the 60s Beurling (see \cite{Be89}) showed that the only $L^1$ function satisfying
\[\iint_{\mathbb{R}\times\mathbb{R}}f(x)\widehat{f}(\xi)e^{|x\xi|}\,dxd\xi<+\infty
 \]
is the identically zero function, property that implies Hardy's uncertainty principle and the generalizations by Cowling and Price, also Morgan's uncertainty principle. This is known in the literature as Beurling's or Beurling--H\"ormander's uncertainty principle, as the proof of this result was published by H\"ormander, \cite{Ho91}.
\item  Sitaram, Sundari and Thangavelu extend Hardy's uncertainty principle to higher dimensions, replacing $x^2$ by $|x|^2$, and to the Heisenberg group, \cite{SST95}. 
 \end{itemize} 

These different extensions, as well as Morgan's uncertainty principle, meet in the work \cite{BDJ03} of Bonami, Demange and Jaming, who extend Cowling--Price and Beurling--H\"ormander theorems to $\mathbb{R}^d$. Also, the authors weakened the decay assumptions so that Hermite functions meet them, very much in the spirit of Theorem \ref{generalhup}. Let us state the general result in its $L^1$ version as originally stated in \cite{BDJ03}, since all the statements in $\mathbb{R}$ or $\mathbb{R}^d$ mentioned so far can be seen as consequences of it:
\begin{thm}\label{bdjbeurling}
    Let $f\in L^2(\mathbb{R}^d)$ and $N\ge0.$ Then
    \begin{equation}\label{conditionbdj}
\iint_{\mathbb{R}^d\times\mathbb{R}^d}\frac{|f(x)||\widehat{f}(\xi)|}{(1+\|x\|+\|\xi\|)^N}e^{|\langle x,\xi\rangle|}\,dxd\xi<+\infty
    \end{equation}
    if and only if $f$ may be written as $f(x)=P(x)e^{-\langle Ax,x\rangle/2}$ where $A$ is a real positive definite symmetric matrix and $P$ is a polynomial of degree at most $\dfrac{N-d}{2}.$
\end{thm}

As this result is a generalization of the previous ones, let us review the proof.

\begin{proof}[Sketch of the proof] The proof is a consequence of three main facts.

 \textbf{Fact 1.} \textit{It is enough to check that \eqref{conditionbdj} implies that $f(x)$ may be written as
\[
f(x)=P(x)e^{-\langle (A+iB)x,x\rangle/2}
\]
where $A$ and $B$ are symmetric matrices and $P$ is a polynomial.}

Indeed, once this is checked, after a change of variables, one may assume that $A$ is the identity matrix, and then
\[
f(x)=P(x)e^{-(\|x\|^2+i\langle Bx,x\rangle)/2}\Longleftrightarrow \widehat{f}(\xi)=Q(\xi)e^{-\langle(I+i B)^{-1}\xi,\xi\rangle/2}
\]
with $P$ and $Q$ polynomials of the same degree that we define as $n$. After this assumption, the integrability condition is seen to imply that $\|x\|^2+\langle(I+B^2)^{-1}\xi,\xi\rangle-2\langle x,\xi\rangle$ is always positive, which is only possible if $B=0.$ To prove the condition on the degree of the polynomial, one takes a point $x_0\in\R^d$ of norm $1$ for which $P_n$ and $Q_n$, the homogeneous terms of maximal degree of $P$ and $Q$ satisfy $P_n(x_0)Q_n(x_0)\ne0$. This allows to prove that in a cone $\Gamma$ of amplitude small enough whose generating line contains $x_0$, $|P_n(x)Q_n(\xi)|\ge c\|x\|^n\|\xi\|^n$. If $x$ and $y$ are large enough, the same inequality holds true for $P$ and $Q$. By \eqref{conditionbdj},
\[
\iint_{\Gamma\times\Gamma} \frac{\|x\|^n\|\xi\|^n}{(1+\|x\|+\|\xi\|)^N}e^{-\|x-\xi\|^2/2}<+\infty.
\]

Another change of variables and Fubini show that there exists $t$ for which
\[
\int_{\Gamma} \frac{\|x\|^n\|x+t\|^n}{(1+\|x\|+\|x+t\|)^N}\,dx<+\infty,
\]
only possible if $2n<N-d$.

 \textbf{Fact 2.} \textit{The function $g$ defined by $\widehat{g}(\xi)=\widehat{f}(\xi)e^{-\|\xi\|^2/2}$ admits an holomorphic extension to $\C^d$ of order $2$ (i.e, there exist constants $c,C>0$ such that $|g(z)|\le Ce^{c\|z\|^2}$ for all $z\in\C^d$). Moreover, there exists a polynomial $R$ such that for all $z\in\C^d,\ g(z)g(iz)=R(z)$.}

The fact that $g$ has an holomorphic extension of order $2$ is a simple consequence of $\widehat{f}\in L^\infty(\R^d)$ thanks to the decay condition. 
 
 The interesting part is to derive the functional equation satisfied by $g$, a nontrivial consequence of Phragm\'en--Lindel\"of's principle and
\begin{equation}\label{polyn}
\int_{\|x\|\le R}\int_{\R^d}|g(x)|\widehat{g}(\xi)|e^{|\langle x,\xi\rangle|}\,dxd\xi<C(1+R)^N,
\end{equation}
which follows from a careful study of the decay condition.

In order to apply Phragm\'en--Lindel\"of's principle, let us take $\alpha\in(0,\pi/2),\ \xi\in\R^d$ and define, for $z\in \C$,
\[
G_\xi^{(\alpha)}(z)=\int_0^{z\xi_1}\cdots \int_0^{z\xi_d}g(e^{-i\alpha }u)g(iu)\,du.
\]

We point out that this is a function of just one complex variable. By Fourier inversion, one checks thanks to \eqref{polyn} that $G_\xi^{(\alpha)}$ has polynomial growth of order $N$ on $e^{i\alpha}\R$ and on $i\R$, with uniform estimates in $\alpha$. By Phragm\'en--Lindel\"of's principle, the same estimate holds in the angular sector
\[
\{r e^{i\beta}: r\ge0,\alpha\le \beta\le \pi/2\},
\]
and, letting $\alpha\to0$, we have an analogous estimate in the first quadrant for $G_\xi^{(0)}$. Working analogously in the other three quadrants, this is an entire function with polynomial growth of order $N$, so a polynomial of degree at most $N$. This allows us to write
\[
G_\xi^{(0)}(z)=a_0(\xi)+\cdots+a_N(\xi)z^N,
\]
being $a_j$ homogeneous polynomials of degree $j$ in $\R^d.$ Furthermore $G_\xi^{(0)}(1)$ is a polynomial on $\R^d$ and admits an entire extension to $\C^d$, that is also a polynomial. By differentiating this extension, we get that $g(z)g(iz)$ is a polynomial.

 \textbf{Fact 3.}  Inequality \textit{\eqref{conditionbdj} implies that $f(x)$ may be written as
\[
f(x)=P(x)e^{-\langle (A+iB)x,x\rangle/2}
\]
where $A$ and $B$ are symmetric matrices and $P$ is a polynomial.}

 The authors show that entire functions $\varphi$ of order 2 of $d$ variables such that on every complex line either they are identically zero or they have at most $N$ zeros are written as $P(z)e^{Q(z)}$ with $P$ a polynomial of at most degree $N$ and $Q(z)$ a polynomial of degree at most 2.

Here again one defines a function of just one variable $\varphi_z(t)=\varphi(zt)$, that is, in this case, an entire function of order $2$ with at most $N$ zeros. By Hadamard's factorization theorem, for every $z\in\C^d$,
\[
\varphi(tz)=(a_0(z)+\cdots+a_N(z)t^N)e^{\alpha(z)t+\beta(z)t^2},
\]
where $\alpha$ and $\beta$ are homogeneous functions of degree 1 and 2 respectively, and $a_j$ is of degree $j$. A careful study of this identity allows us to check that $\alpha,\beta$ and $a_j$ are all polynomials.

Finally, an application of this result to the previous function $g$ implies that $g(z)=P(z) e^{Q(z)}$ with $P$ a polynomial of degree at most $N$ and $Q$ a polynomial of degree at most 2. Further, $g(z)g(iz)$ being a polynomial implies that $Q(z)+Q(iz)=0$, so $Q$ is a homogeneous polynomial of degree 2. Hence, $g$ is written in the desired form, and so is $f$.
\end{proof}

As explained by the authors, their arguments are inspired by the proof of H\"ormander (case $N=0,\ d=1$). Also, $L^p$ versions of this theorem also hold true. Interestingly, this general theorem allows for the following modification of Theorem \ref{generalhup}, also discussed in \cite{BDJ03}:

\begin{thm}\label{bdjhardy}
    Let $f\in L^2(\mathbb{R}^d)$ be such that
    \[
|f(x)|\le C(1+\|x\|)^N e^{-\langle Ax,x\rangle/2}\ \ \text{and}\ \ |\widehat{f}(\xi)|\le C(1+\|\xi\|)^N e^{-\langle B\xi,\xi\rangle/2},
    \]
    where $A$ and $B$ are two real symmetric positive definite matrices. If $B-A^{-1}$ is semi-positive definite and non-zero, then $f$ is null. If $B=A^{-1}$, then $f(x)=P(x)e^{-\langle Ax,x\rangle/2}$ where $P$ is a polynomial of degree at most $N$. Else, there is a dense subspace of functions satisfying these estimates.
\end{thm}

We want to emphasize that the main tools of the proof discussed above, and also of the proofs of all the statements discussed so far, are based on complex analytic tools.

\section{Hardy's uncertainty principle and real variable arguments}

A natural question arises: Is there a real variable proof of Hardy's uncertainty principle? The search for this proof is what has motivated several relevant results concerning the impact of Hardy's paper. Independently and almost simultaneously we find two related results.

On the one hand, Tao reviews in his blog \cite{Ta10} Hardy's uncertainty principle and is curious about a real variable proof of it, at least in the one dimensional case. He exploits estimates for the derivatives of the Fourier transform of a function with Gaussian decay, together with doubling bounds for polynomials, to provide a weaker version of Hardy's uncertainty principle. In the context of Remark \ref{hardyab} he shows the existence of a constant $C_0$ such that $ab>C_0$ implies that the function is trivial. His proof does not quantify the constant $C_0$ and he suggests that most likely the best constant that his proof allows is not $1/4$, the sharp value.

On the other hand, Kenig, Ponce, and the second author in \cite{KPV03}, and after this article, together with Escauriaza in \cite{EKPV06,EKPV07}, start to study unique continuation properties for solutions to dispersive equations, mainly Korteweg-de-Vries and Schr\"odinger equations. We can, roughly speaking, summarize their statements in the following form:

What is the maximum decay that a nontrivial $u$, a strong solution of a certain PDE $Pu=0$, can admit?

It turns out that the type of results they look for resembles Hardy's uncertainty principle. Independently and, again, almost simultaneously, Chanillo makes this connection clear in \cite{Ch07}, showing an equivalence between Hardy's uncertainty principle and a uniqueness result for Schr\"odinger evolutions in connected complex semi-simple Lie groups. In order to discuss Chanillo's result, let us introduce some notation, borrowed from \cite{Ch07}.
Let $G$ denote a complex, connected, semi-simple Lie group, $K$ a fixed maximal compact subgroup of $G$, and $\mathcal{G}$ and $\mathcal{K}$ their Lie algebras respectively. Let $B$ denote the Cartan--Killing form on $\mathcal{G}$, and the Cartan decomposition is given by $\mathcal{G}=\mathcal{K}\oplus\mathcal{P}$. The restriction of $B$ to $\mathcal{P}\times\mathcal{P}$ is strictly positive definite and defines a norm, $|H|^2=B(H,H)$. Let $\mathcal{A}$ be a fixed maximal abelian subspace of $\mathcal{P}$. Let $\Sigma$ denote the set of non-zero roots corresponding to the pair $(\mathcal{G},\mathcal{A})$, and $\Sigma_+$ the set of positive roots $\alpha$ for some ordering. Let $A$ denote the analytic subgroup with Lie algebra $\mathcal{A}$ and $W$ the Weyl group associated to $\Sigma$.
The Haar measure on $A$ can then be written by a formula of Haris--Chandra as
\[
dx=\left(\sum_{s\in W}(\det s)e^{s\rho(H)}\right)^2\,dH=(\phi(H))^2\,dH,\ H\in\mathcal{A}.
\]
Here $\rho=\frac12\sum_{\Sigma_+}\alpha$. The spectral variables are elements of $\mathcal{A}^*$ and will be denote by $\lambda$. On $G/K$ we have a $G$ invariant Riemannian metric and using this metric we can form a Laplace--Beltrami operator $\Delta$. Moreover, there is a unique element $H_\lambda\in\mathcal{A}$ such that
\[
\lambda(H)=B(H_\lambda,H),\ H\in\mathcal{A}.
\]
Finally, we say that a function $f$ on $G$ is $K$ bi-invariant if and only if $f(k_1 a k_2)=f(a),\ \forall k_1,k_2\in K,\ \forall a\in A.$
Then, Chanillo's result is the following:
\begin{thm}\label{thmchanillo}
Let $G$ denote a connected, complex, semi-simple Lie group. Let $f$ be a bi-invariant function such that
\[
|f(H)\le A e^{-a|H|^2}.
\]
Consider the initial value problem,
\[
-i u_t=\Delta u,\ u_{|_{t=0}}=f.
\]
Then, $u$ is also bi-invariant. Furthermore, if
\[
|u(H,t_0)|\le B e^{-b|H|^2},
\]
and $16abt_0^2>1,$ then necessarily $u\equiv0$ for all $t\ge0$.
\end{thm}
The main part of the proof relies on getting the formula
\[
u(H_1,t)\phi(H_1)=c_l |W|^2 \frac{e^{-i(t|\rho|^2-\frac{|H_1|^2}{4t})}}{t^{l/2}}\widehat{R}(\frac{H_1}{2t}),
\]
where $l=\dim \mathcal{A}$, $R(H)=e^{i\frac{|H|^2}{4t}}f(H)\phi(H)$, and the hat stands for the Fourier transform in $\R^l$.
Noting that $|\phi(H)|\le c e^{c|H|}$, using the decay condition at $t=t_0$ and the definition of $R$ one gets
\[
|\widehat{R}(H)|\le e^{-4b't_0^2|H|^2},\ |R(H)|\le c e^{-a|H|^2},
\]
where $b'<b$ and $a'<a$. By Hardy's uncertainty principle, we conclude that $R\equiv0$ and hence $f\equiv0$.
We point out that the formula for the solution of Schr\"odinger evolution is a generalization of the well known formula
\begin{equation}\label{freeschr}
u(x,t)=c_d \frac{e^{i|x|^2/(4t)}}{t^{d/2}}\left(e^{i|y|^2/(4t)}f\right)\widehat{}\left(\frac{x}{2t}\right).
\end{equation}
This formula may also be used to generalize Hardy's uncertainty principle in the directions discussed in the previous section, and also to show the end-point case: If $16abt_0^2=1$, then $u(x,0)$ is a multiple of $e^{-a|x|^2-i|x|^2/(4t_0)}$. Identity \eqref{freeschr} gives dynamical interpretations of uncertainty principles, which can be seen as uniqueness results of solutions to the free Schr\"odinger equation. In this case, having a real calculus proof of Hardy's uncertainty principle would most likely allow one to extend this dynamical version to more general (perturbed) problems.

%\blue{It turns out that the type of results they look for resembles Hardy's uncertainty principle. Independently and, again, almost simultaneously, Chanillo makes this connection clear in \cite{Ch07}, showing that Hardy's uncertainty principle implies:
%\begin{thm}\label{thmchanillo}
%    Let us consider the initial value problem for $u(x,t)$,
%    \[-iu_t=\Delta u,\ u|_{t=0}=f(x).\]
%    Assume that $|f(x)|\le Ae^{-a|x|^2}$, and $|u(x,t_0)|\le Be^{-%b|x|^2}$. If $16abt_0^2>1$, then $u(x,t)\equiv0$, for all $t\ge0.$ 
%\end{thm}
%This result, stated in \cite{Ch07} for connected complex semi-simple Lie groups, is a consequence of the formula
%\begin{equation}\label{freeschr}
%u(x,t)=c_d \frac{e^{i|x|^2/(4t)}}{t^{d/2}}\left(e^{i|y|^2/(4t)}f\right)\widehat{}\left(\frac{x}{2t}\right),
%\end{equation}
%which can be applied not only to prove Theorem \ref{thmchanillo}, but also to generalize Hardy's uncertainty principle in the directions discussed in the previous section, and to show the end-point case: If $16abt_0^2=1$, then $u(x,0)$ is a multiple of $e^{-a|x|^2-i|x|^2/(4t_0)}$. Identity \eqref{freeschr} gives dynamical interpretations of uncertainty principles, which can be seen as uniqueness results of solutions to the free Schr\"odinger equation. In this case, having a real calculus proof of Hardy's uncertainty principle would most likely allow to extend this dynamical versions to more general (perturbed) problems.} 

The work of Escauriaza, Kenig, Ponce and the second author can be undertstood in the opposite direction. They provide, in a PDE setting, real variable proofs of dynamical versions of Hardy's (or Morgan's) uncertainty principle. The main motivation behind this research is the fact that adding a general potential to the Schr\"odinger equation does not allow a direct argument to be used as \eqref{freeschr} to prove an analogous version of Theorem \ref{thmchanillo}, so a new argument is required. Their methods allow for a complete proof of Hardy's uncertainty principle, as stated in Remark \ref{hardyab} and adapt to a large variety of settings, as we will describe in the next section. 

The strategy of proof combines two facts. On the one hand, a Carleman estimate, which is an $L^2$-weighted inequality of the form
\[
\|e^{\varphi}f\|_{L^2(\mathbb{R}^{d+1})} \le C_\varphi\|e^{\varphi}(i\partial_t+\Delta)f\|_{L^2(\mathbb{R}^{d+1})},
\]
for any \textit{reasonable} function $f$ and an appropriate weight function $\varphi$. On the other hand, convexity properties of solutions to the Schr\"odinger equation with fast decay at two times to prove, roughly speaking, how the fast decay not only is preserved to any time in between, but gets better in the interior of the interval limited by the two times when one has fast decay.

These two properties allow the authors to give a lower and an upper bound for the $L^2$-norm of a nonzero solution that leads to a contradiction if the decay hypothesis at two times are appropriately chosen. More precisely, they give in \cite{EKPV10} an alternative proof of Theorem \ref{thmchanillo} when $G=\R^d$ stating the decay hypothesis in an $L^2$ sense. This proof allows to consider more general equations, $-iu_t(x,t)=\Delta u(x,t)+V(x,t)u(x,t)$, under size constraints in the potential $V(x,t)$.

We can also understand this uniqueness result as a non-stationary version of Landis' conjecture. In 1988, \cite{KL88}, Kondrat'ev and Landis conjectured that a nontrivial solution to $(-\Delta +V(x))u=0$ cannot satisfy the estimate
\[
|u(x)|\le C e^{-(\|V\|^{1/2}_{L^\infty(\R^d)}+\epsilon)|x|}
\]
with $C$ and $\epsilon$ positive constants. The conjecture was disproved by Meshkov, \cite{Me89}, who constructed a counterexample. However, the conjecture is still open in full generality if we impose the potential to be real valued, despite the recent progress by Logunov, Malinnikova, Nadirashvili and Nazarov \cite{LMNN20} in dimension two. For more information about Landis' conjecture, we refer the reader to the recent survey \cite{FBRS24}.

Going back to the dynamical uncertainty principle, the real variable proof is a consequence of a series theorems, each one improving the result of the previous theorem, thanks to a better understanding of solutions to the Schr\"odinger equation with fast decay at two times. Eventually, not just a full proof of the result stated in Remark \ref{hardyab} is obtained but also what could be called a particular quantitative version of Hardy's uncertainty principle, very much in the spirit of the Heisenberg uncertainty principle

A fundamental part of the proof is to reduce the proof of the Carleman estimate to a logarithmic convexity argument that goes back to Agmon and Nirenberg \cite{Ag66}.

Consider, as we did in Section 2  two operators, $\mathcal S $ symmetric and $\mathcal A $ skewsymmetric. We look at solutions of the abstract evolution equation:
$$\partial_t v=(\mathcal S+\mathcal A)v.$$
For simplification we will assume that $\mathcal S $  and $\mathcal A $ are independent of time. Our interest is the function
$$H(t)=\langle v,v\rangle.$$
Formally, we get
$$\begin{aligned}
\dot H(t)&=\langle v_t,v\rangle+\langle v,v_t\rangle\\
&=\langle(S+\mathcal A)v,v\rangle+\langle v, (S+\mathcal A)v\rangle\\
&=2 \langle Sv,v\rangle.\end{aligned}.$$
Similarly
$$\begin{aligned}\ddot H(t)&=2\langle S v_t,v\rangle+2 \langle S v,v_t\rangle\\
&=2\langle(S+\mathcal A)v,S v\rangle+2\langle S v, (\mathcal A+S)v\rangle\\
&=4 \langle Sv,Sv\rangle+2\langle (S\mathcal A-\mathcal A S)v,v\rangle.
\end{aligned}$$
As a consequence,
$$\begin{array}{rcl}
\big(\log H(t)\big )^{\cdot\cdot}&=&\left(\displaystyle\frac{\dot H}{H}\right)^{\cdot}=\displaystyle\frac{\ddot H H-\dot H^2}{H^2}\\
&=&\displaystyle\frac{1}{\langle v, v\rangle}\left\{4\langle \mathcal S v,\mathcal S v\rangle \langle v,v\rangle-4\langle \mathcal S v,v\rangle^2+2\left\langle(\mathcal S\mathcal A-\mathcal A\mathcal S)v,v\right\rangle\right\}\\
&\ge& 2\dfrac{\left\langle (\mathcal S \mathcal A-\mathcal A \mathcal S)v,v\right\rangle}{\left\langle v,v\right\rangle}.
\end{array}$$
Hence, if 
$$2\left\langle (\mathcal S\mathcal A-\mathcal A \mathcal S)v,v\right\rangle\ge 0\,$$ 
then
\begin{equation*}
H(t)\le H(0)^{1-t} H(1)^t.
\end{equation*}
More generally, if 
\begin{equation*}
2\left\langle (\mathcal S\mathcal A-\mathcal A \mathcal S)v,v\right\rangle\ge \psi (t) \langle v,v\rangle,
\end{equation*}
then
\begin{equation*}
H(t) e^{-B(t)}\le H(0)^{1-t} H(1)^t,
\end{equation*}
with 
\begin{equation*}
\ddot B=\Psi(t),\quad B(0)=B(1)=0.
\end{equation*}

One relevant example of the above type of argument at the heart of the results in \cite{EKPV06, EKPV08,EKPV10} is the following. Take $u$ a solution of
$$\partial_t u=i \Delta u$$
and define
$$e^{\frac{|x|^2}{2}}u=v$$
and
$$H(t)=\|v(t)\|_{L^2}^2=\langle v(t),v(t)\rangle.$$
Then formally,
$$\begin{array}{rcl}
\partial_t v&=&\left(e^{\frac{|x|^2}{2}}i \Delta e^{-\frac{|x|^2}{2}} \right)v\\
&=&i \displaystyle\sum\limits_j e^{\frac{|x|^2}{2}}\partial^2_je^{-\frac{|x|^2}{2}}v\\
&=&i \displaystyle\sum\limits_j e^{\frac{|x|^2}{2}}\partial_je^{-\frac{|x|^2}{2}} e^{\frac{|x|^2}{2}}\partial_je^{-\frac{|x|^2}{2}} v\\
&=& i \displaystyle\sum\limits_j (x_j-\partial_j) (x_j-\partial_j) v\\
&=& i\big(|x|^2-2x\cdot\nabla-d+\Delta\big) v.
\end{array}$$
This implies that $v$ satisfies
\begin{equation*}
\partial_t v=(\mathcal S+\mathcal A)v\quad\text{with}\quad \mathcal S v=-i(2x\cdot\nabla+d)v\quad\text{and}\quad \mathcal  Av=i(\Delta+|x|^2).
\end{equation*}
A simple computation gives
\begin{equation*}
\mathcal S \mathcal  A-\mathcal  A \mathcal S=-4\Delta+4|x|^2,
\end{equation*}
and from the Heisenberg uncertainty principle
\begin{equation*}
\big\langle \mathcal S \mathcal  A-\mathcal  A \mathcal S v, v\big\rangle=4\big(\|\nabla v\|^2+\|xv\|^2\big)\ge 4\langle v,v\rangle.
\end{equation*}
Hence as we proved above
\begin{equation*}
\big(\log H(t)\big)^{\cdot\cdot}\ge 8
\end{equation*}
and
\begin{equation*}
H(t)\le H(0)^{1-t} H(1)^t.
\end{equation*}
As a consequence, if $H(0)$ and $H(1)$ are finite then at a formal level $H(t)$ is also finite and  $u$ has Gaussian decay for $0<t<1$. Observe that if $u_0$ has Gaussian decay, then $e^{it\Delta}u_0$ does not necessarily have it for $t>0$ as can be easily checked taking $u_0=({\textrm{sig}\,}x)e^{-|x|^2}$, with ${\textrm{sig}\,}x$ denoting the signum function. To go from the formal level to a rigorous one is a delicate question and is an important part of the arguments in \cite{EKPV06, EKPV08, EKPV10}.

Formally speaking, the argument presented above allows to, in a first step, deal with the case where a decay of the form $e^{\gamma|\cdot|^\alpha}u(\cdot,s_i)\in L^2(\mathbb{R}^d)$, with $\gamma>0$ and $\alpha>2$, $i\in\{1,2\}$ is assumed or, to give a closer relation to Hardy's uncertainty principle, one can replace this assumption by $e^{\gamma|\cdot|^2}u(\cdot,s_i)\in L^2(\mathbb{R}^d)$ for a certain large constant $\gamma$, which cannot be quantified. This is the main result in \cite{EKPV06}.

In a second step, a better use of Carleman estimates and convexity estimates provides in \cite{EKPV08} a quantification of the latter constant $\gamma$, extending Theorem \ref{thmchanillo} to the perturbative setting if $8abt_0^2>1$. 

The last step consists of reducing the gap of the previous result towards the condition coming from the free case. This reduction comes from an iterative process, explained in \cite{EKPV10}, where at each step the decay assumption implying that the solution is identically zero gets closer to $16abt_0^2>1$. It is seen that this iterative process is infinite and at the limit one gets the desired restriction on $a$ and $b$. 

This result is sharp, in the sense that in the case $16abt_0^2=1$ one can construct a potential and a nontrivial solution satisfying the corresponding decay condition. Also, through $\eqref{freeschr}$, neglecting the perturbations, the authors provide a real variable proof of part of Hardy's uncertainty principle, as their result does not recover either the end-point case of the $L^2$ version ($16abt_0^2=1$ implies $u\equiv0$) or the end-point case of the $L^\infty$ version discussed above. It is also relevant to note that Hardy's uncertainty principle is a qualitative uncertainty principle, in the sense that the identically zero function is the only one satisfying certain decay conditions. Chanillo's result in Theorem \ref{thmchanillo} also exhibits a qualitative behavior. However, the real variable approach we have briefly described here gives a quantitative result. This is because the same iterative process can be run in the case $16abt_0^2<1$. In this case, the iterative process is still infinite, but at the limit one gets that a nontrivial solution with Gaussian decay at times $t=0$ and $t=t_0$ satisfies a precise Gaussian decay in the interval $(0,t_0)$. For a more detailed description of this quantitative statement we refer to \cite{EKPV10}.

Regarding the end-point case, Escauriaza, Kenig, Ponce and the second author, together with Cowling include modifications of the previous arguments in order to revisit Hardy's uncertainty principle in \cite{CEKPV10}, providing a full real variable proof of Theorem \ref{classicalhup}, up to the end-point case. We point out that the arguments performed in the PDE setting adapt very well to the $L^2$ version of Hardy's uncertainty principle, in contrast with the $L^\infty$ formulation of Theorem \ref{classicalhup}. A nice trick in \cite{CEKPV10} concludes the end-point case of the $L^\infty$ version from the one of the $L^2$ version:

In one dimension, if $ab=1/4$ and, 
\[
\|e^{ax^2}f\|_{L^\infty(\R)}+\|e^{b\xi^2}\widehat{f}\|_{L^\infty(\R)}<+\infty,
\]
then $f$ is smooth, as well as $g(x)=\frac{f(x)-f(0)e^{-ax^2}}{x}$, and its Fourier transform, which is, up to a multiple, the convolution of the sign function with $(f(\cdot)-f(0)e^{-a(\cdot)^2})^{\widehat{}}$. Careful analysis of this convolution implies that $g$ satisfies the end-point case of the $L^2(\R)$ version of the Hardy uncertainty principle, so $g$ is identically zero and, therefore, $f(x)=f(0)e^{-ax^2}$. The higher-dimensional case follows easily from the one-dimensional case.

Some years later a new one-dimensional real variable proof was proposed by Pauwels and de Gosson in \cite{PdG14}. In this case, it does not cover the full range of values $a$ and $b$, as it requires $L^\infty(\R)$ estimates under the constraint $ab\ge 1$. Their proof connects Hardy's uncertainty principle to the theory of band limited functions, as the key idea is to understand the asymptotic behavior of the largest eigenvalue of the sum of the time and band limited operators as both the time limit and the band limit approach infinity.

Denoting the time limit operator and the band limit operator respectively as $\chi_{(-\tau,\tau)}$ and $S_\omega$, they relate the eigenvalues of the operator $\chi_{(-\tau,\tau)}+S_\omega$ to the eigenvalues of the operator
\[
Pf(x)=\int_{-1}^1\frac{\sin c(x-y)}{\pi(x-y)}f(y)\,dy,
\]
whose eigenfunctions are the so called prolate spheroidal wave functions, well studied in the 60s by Landau, Pollak, and Slepian.

\section{Hardy's uncertainty principle in different settings}

We would like to finish this short survey by pointing out some recent extensions of the Hardy's uncertainty principle to different settings beyond $\R^d$, as all that we have previously discussed refers to $\R^d$. The road-map of the real variable proof in \cite{EKPV10} has been seen to be very robust and applicable to a large variety of settings. 

Still in the dynamical context related to Schr\"odinger evolutions in $\R^d$, Dong and Staubach included a drift term in \cite{DS07}, in order to study unique continuation properties of solutions to
\[
u_t=i(\Delta u+Vu +B\cdot\nabla u),
\]
showing that under size conditions on the potentials $V,B$ a nontrivial $C([0,1],H^1(\R^n))$-solution to the previous equation cannot satisfy $u(\cdot,0),u(\cdot,1)\in H^1(e^{a|x|^2})$ for $a$ large enough. A quantization of the sharp value of $a$ is missing, as well as a relaxation of the hypothesis from $H^1$ to $L^2$. While, up to our knowledge, there has not been any progress in a context with such a generality, the somewhat related case of the magnetic Schr\"odinger equation
\[
u_t=i((\nabla-iA(x))^2+V)u,\ A:\R^n\to\R^n
\]
has been studied in great detail in \cite{BFGRV13, CF15}, providing a sharp result. Also, as an application of this case, harmonic oscillators and uniform magnetic fields are considered in \cite{CF17}. We also point out the work of Knutsen, extending the method to Schr\"odinger equations with real quadratic Hamiltonians in \cite{Kn23}, by using an extension of Hardy's uncertainty principle to the Wigner distribution.

Escauriaza and collaborators extended their method in \cite{EKPV16} to provide a dynamical interpretation of Hardy's uncertainty principle seen as a uniqueness result for solutions to the heat equation. In this parabolic setting, starting from a suggestion by Zuazua and the survey \cite{Ja06} by Jaming, they see that there is no nontrivial solution to a perturbed heat equation with $L^2$ initial datum that has fast decay (also related to a Gaussian) at a posterior time. Similar results for more generic parabolic problems are studied in \cite{Ng10}, and, also, under the presence of stochastic perturbations, \cite{FBZ20}.

The research of the first author of this short note in collaboration with different researchers has been focused on finding analogous results in different settings. We refer to the survey \cite{FBM21} and the references therein, where dynamical discrete versions of Hardy's uncertainty principle in $\Z^d$, in homogeneous trees and on quantum regular trees can be found. In all these scenarios, sharp results are yet to be proved.

In the discrete setting, it is also interesting to point out the works \cite{JLMP18} by Jaming, Lyubarskii, Malinnikova and Perfekt and \cite{LM18} by Lyubarskii and Malinnikova in $\Z^d$, also the works by Álvarez-Romero and Teschl in different graphs and for different evolutions \cite{AV18,AVT17a,AVT17b}. Using complex analytic tools we have a sharp uniqueness result for solutions to the free discrete Schr\"odinger with fast decay at two times, given the decay in terms of the modified Bessel function. 

Putting aside the dynamical interpretations of Hardy's uncertainty principle and going back to the original statement, de Gosson gives in \cite{dG17} geometrical interpretations of the principle. In the one-dimensional case these interpretations may be rather intuitive, as we can rewrite Remark \ref{hardyab} in the following geometric way:
\textit{There exist nontrivial functions $f$ satisfying $f=O(e^{-ax^2})$ and $\widehat{f}=O(e^{-b\xi^2})$} if the area of the ellipse $\Omega_{a,b}=\{(x,\xi)\in\R^2: ax^2+b\xi^2 \le 1\}$ is bigger than $2\pi$.
The question under study is to find suitable analogs in the higher dimensional case.

Another line of research developed since the millennium is the extension of Hardy's uncertainty principle to other operators besides the Fourier transform. We have mentioned above the role of the Wigner distribution in \cite{Kn23}, based on the work by Gr\"ochenig and Zimmermann \cite{GZ01}, who, under the conviction that every uncertainty principle about $f$ and $\widehat{f}$ can be translated into an equality about the short-time Fourier transform or the Wigner distribution, extend Hardy's uncertainty principle to these two transformations. In the last few years new perspectives of uncertainty principles have been considered, with special interest in metaplectic operators as these operators have applications in modern aspects of Analysis and PDEs, such as time-frequency analysis or quantum harmonic analysis. Extensions of several uncertainty principles are covered in \cite{DdGP24}. We also point out to the very recent paper \cite{CGM24}, connecting these results to Schr\"odinger evolutions in the same spirit as \cite{Kn23}.

Finally, we have already mentioned that Hardy's uncertainty principle is extended in \cite{SST95} to the Heisenberg group; the latter paper tackles the question of extending the principle to more general Lie groups, and it is the first of a large amount of extensions, as we can see in \cite{BK08,BT10,KK01,PT09,Th01,Th02} replacing the Gaussian by the heat kernel, and, summarizing the main results under the sentence \textit{only $0$ has faster decay than the heat kernel}. We have already seen that Chanillo in \cite{Ch07} extends the principle to Lie groups by relating the problem to a uniqueness problem for the free Schr\"odinger evolution with fast decay at two times.
%\blue{The above cited work of Chanillo \cite{Ch07} not only characterizes Hardy's uncertainty principle as a uniqueness result for the free Schr\"odinger evolution with fast decay at two times, but uses this characterization to extend the principle to Lie groups.} 
In this dynamical direction, more generalizations to the Heisenberg group and also to $H$-type groups are studied in \cite{BTD13,FBJPE21,JG25,LM14}.

\end{document}